*Optimal scheduling of energy and mass flows based on networked multi-carrier hubs formulation: a general framework*
*Mohamed Tahar Mabrouk[a], Mahdi Majidniya[b] and Bruno Lacarrière[c]*

[a] IMT Atlantique, Department of Energy Systems and Environment, GEPEA, UMR CNRS, 6144, F-44307, Nantes, France, mohamed-tahar.mabrouk@imt-atlantique.fr, CA
[b] IMT Atlantique, Department of Energy Systems and Environment, GEPEA, UMR CNRS, 6144, F-44307, Nantes, France, mahdi.majidniya@imt-atlantique.fr
[c] IMT Atlantique, Department of Energy Systems and Environment, GEPEA, UMR CNRS, 6144, F-44307, Nantes, France, bruno.lacarriere@imt-atlantique.fr

**Abstract:**
Due to increased energy demand and environmental concerns such as greenhouse gas emissions and natural resources depletion, optimizing energy and raw materials usage has recently drawn much attention. Achieving more synergy between different energy sectors and manufacturing processes could lead to substantial improvements. Energy hubs are already well-known solutions for studying multi-carrier energy systems. In the present study, a multi-carrier hub is defined as a geographic area where different processes take place to convert energy and material flows possibly consumed locally, stored, or exported. Hub boundaries correspond to an area small enough to neglect losses when energy and materials flows are exchanged between the processes. A general formulation is introduced to model such a hub. This formalism encompasses all the possible configurations without limitations or preconceptions regarding used carriers or processes' architecture. For this purpose, a new hub representation is proposed that allows optimization while considering all possible arrangements (parallel, serial, or a combination of them). It is implemented in a framework that allows the creation of several multi-carrier hubs and connecting some of them to exchange flows through dedicated networks. The formalism is based on linear programming (LP). The features of the developed framework are illustrated by a case study consisting of a network made of three hubs and considering heat, electricity, water, and hydrogen flows.The study aims to optimally schedule processes usage and materials and energy storage to minimize the cost of imported resources. The study is done in a short-term period. It shows the benefits of synergy between the different processes in the hub.

# 1. Introduction

The increasing demand for energy and raw materials requires corresponding improvements in material and energy efficiency in major manufacturing and energy industries to prevent depleting natural resources. Although the needs of the residential and transportation sectors are different, these challenges are similar. Individual process improvements are reaching their limits. Additionally, improving efficiency within a single sector or a small geographical area is also severely limited. Today, finding ways to realize more synergy between different sectors and/or locations where energy and materials are converted, consumed, and stored in a coordinated manner is a topic that is attracting much attention. This is known in some studies as "Urban-Industrial symbiosis"which expends the "Industrial symbiosis"concept and was proven to be beneficial [1] although complex to implement.

## 1.1. Literature review

Optimization models can be powerful tools for finding these new ways of connecting different processes in different sectors and/or geographic locations. Many existing decision-support tools use optimization models and solvers to help users plan and commit distributed production systems. This is the case in the energy field, where plenty of decision-support tools exist. These tools are used for planning and commitment of Multi-Carrier Distributed Energy systems (MC-DES). Many of them use optimization models to compare different technological alternatives and size subsystems. Some of these software tools are proprietary such as DER-CAM and Crystal Energy Planner. The first has been developed and maintained by the Lawrence Berkeley National Laboratory (LBNL), US, since 2000 [2] and it exists in two versions: 1) investment and planning and 2) operation. The second is a commercial software and proposes a techno-economic optimization of energy systems for different time horizons [3] In addition, Open-source tools exist for all problem scales: OSeMOSYS [4] and Energy Scope TD [5] are suitable for large-scale problems ranging from regional to national scale. We can find OEMOF [6] and RIVUS [7] at the district scale. The first one uses a graph approach to model connected multi-carrier technologies with the possibility of using external libraries to increase modularity. The second one focuses on finding the cost-minimal distribution network to satisfy demands for different energy carriers using a graph approach and with the possibility to consume energy along the edges and convert energy carriers in the vertices. Some other tools focus on smaller-scale problems from the building (or factory) to the district scale, such as URBS[8] and its upgraded derivative FICUS [9] and OMEGAlpes [10]. The latter includes an exergy module and an actor package that enables considering social and multi-actor constraints. However, it is limited to electricity, heat, and gas carriers. All these software tools use Mixed-Integer Linear programming (MILP), which is a powerful tool to solve large and complex energy integration problems. However, this method has several limitations that led researchers to develop nonlinear models using deterministic approaches [11], stochastic approaches [12], or bi-level optimization techniques [13,14].

On the other hand, plenty of different formalisms are used in these tools, such as graph formalisms (OEMOF and RIVUS). OMEGAlpes uses blocks representing mono-energy units and conversion units connected to each other using energy flow connections and nodes. URBS and FICUS are based on the same principle. They use several model entities such as processes, commodities, storage, and intermittent supply. The main conceptual difference between them is that URBS allows flow transmission between different sites, whereas FICUS considers one factory which is defined as a small microgrid with multiple demands and energy-conversion technologies. This concept is very close to the Energy Hub (EH), widely used for the optimal management of multiple energy carrier networks, such as electricity, thermal power, and gas. It was first proposed by Favre-Perrot *et al.* in [15] who defined EH as an interface between producers, consumers, and transmission networks that connect the hubs together. The concept was refined and extended by Gidel *et al.* [16] who defined the EH as a unit where energy carriers are converted, consumed, or stored. They also introduced the conversion matrix as a mathematical formulation of the hub. Since then, the EH concept has been adopted and developed by many researchers. Ayele *et al.* [17] extended the concept by adding internal generation sources and detailed transmission models for heat and electricity networks. More insight about recent advances in the EH concept can be found in [18] and [19]

Recent studies attempted to extend the energy hub concept to what is called "Industrial Hub (IH) by including production scheduling along with energy management [20]. Kantor *et al.* [21] proposed a MILP formulation for optimizing material and energy carriers. In their formulation, the authors define two types of utilities: 1) process units and 2) utility units. The processes belong to different clusters with the possibility to exchange mass flows between clusters for some carriers. These processes exchange flows with the utilities in both directions through 4 clusters: Resources, Networks, Markets, and waste.

## 1.2. Contribution

In this paper, a novel formalism is introduced for modeling Multi-carrier Hubs in which different energy and material carriers are converted, consumed, or stored. The proposed formalism has a degree of generality such that all possible configurations could be modeled in terms of interconnections inside the hubs, between them, or with the outside of the studied network. The hubs can be flexibly connected through transmission networks. The novelty of the present formulation lies in the proposed hub architecture using general balance equations that allows optimization without any preconceptions in terms of connections between the components while providing all possible arrangements (parallel, serial, or a combination of them).

This formalism is implemented in a framework that allows automatic generation of optimization models of Multi-Carrier Hub Networks by simply defining hubs in which processes, loads, storage units, and different types of exchange ports could be added any way the user wishes. The formalism is illustrated by an optimal scheduling problem of three clusters that include several energy conversion and storage systems and a desalination system as an example of material processing.

## 1.3. Paper organization

The paper is organized into five sections. In the first section, the context and the contribution are introduced. In the second section, the theoretical formulation of MCHN is detailed. In section three, the case study is described. In section four, the results of the optimal scheduling of the three hubs network are presented and discussed. The concluding remarks are given in section five.

# 2. Multi-Carrier Hub Network (MCHN): general formulation

## 2.1. Hub description (MCH)

A Multi-Carrier Hub (MCH) is understood here as a system that may convert, consume, or store multiple mass and/or energy flows. Its geographical coverage is small enough to neglect internal distribution losses. The general structure of the hub proposed in this paper is presented in Fig. 1.

This structure is built around the Hub balances bloc, which guarantees the flow balance for each carrier at the hub's level. Indeed, for each carrier, the flow produced by the processes added to the flow discharged from the storage and the flow extracted from the network is either consumed by the corresponding load, stored in the storage devices, exported, injected to the network, or returned to the processes' bloc. The latter gathers all the conversion processes in the hub. A process converts one or more inlet carriers to one or more outlet carriers. Inlet flows going to the processes originate either from outside the hub through the import port or returned from the output balances bloc to the processes. Exporting (resp. importing) is the process of selling (resp. purchasing) the carrier to (from) outside the boundaries of the investigated network, in contrast injecting (resp. extracting) to (from) the network, which is exchanging flow with other hubs within the studied network.

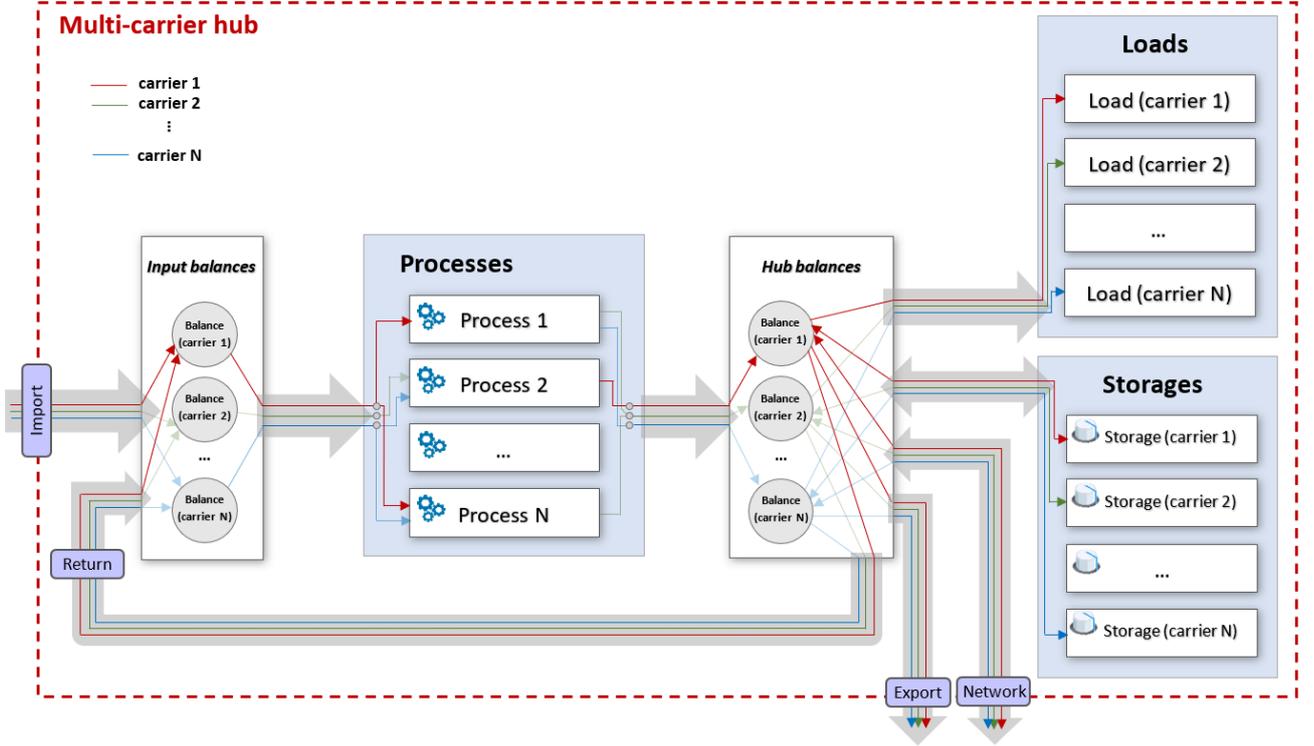

*Fig. 1. The general structure of a Multi-Carrier Hub*

## 2.2. Mathematical formulation

The MCHN is modeled as a constrained optimization problem. The objective function could be adapted to different use cases (optimal operation, optimal sizing of the systems, maximal integration of renewables, etc.) and different objectives (economic, energetic, environmental, etc.). Constraints are used to model the structure and the operation of the hubs, the operation of the conversion processes and the storage systems, and the flow exchange between the hubs through the networks.

It is assumed that the sizes of the conversion processes and the storage units are known. The objective function is the operational cost during the studied period. Also, the generic mathematical formulation proposed here is kept linear. However, it can be easily extended to a nonlinear formulation without changing the problem structure defined in Fig. 1.

Carrier's flows of different types are manipulated in the model (energy and mass flows). For homogeneity purposes, a flow $\mathcal{F}_c$ of a carrier $c$ is defined in this paper as the power flow for energy carriers, and the exergy flow for materials is defined as the product of the mass flow rate $\dot{m}_c$, and a reference chemical exergy $E_{ex}$ as defined in [22] :

$$\mathcal{F}_c = \dot{m}_c \, E_{ex,c} \quad (1)$$

### 2.2.1. Objective function

The objective function is defined as the cost of imported energy and material carriers' flows reduced by the income earned by exporting via the export port:

$$min \left( \sum_{c \in N_h^{imp}} \sum_{h \in N_h} \sum_{t \in N_t} c_{c,h,t} \, \mathcal{F}_{c,h,t}^{imp} \, \Delta t_t - \sum_{c \in N_h^{exp}} \sum_{h \in N_h} \sum_{t \in N_t} p_{c,h,t} \, \mathcal{F}_{c,h,t}^{exp} \, \Delta t_t \right) \quad (2)$$

Where $N_h$ is the set of the hubs in the study, $N_t$ is the set of time steps and $N_h^{imp}$ and $N_h^{exp}$ are the sets of imported and exported carriers in a hub $h$.

$\mathcal{F}_{c,h,t}^{imp}$ and $\mathcal{F}_{c,h,t}^{exp}$ are respectively imported and exported flows of the carrier $c$ in the hub $h$ at the time step $t$. $c_{c,h,t}$ and $p_{c,h,t}$ are the cost (resp. the price) of imported (resp. exported) flows of the carrier c in the hub h at the time step t.

### 2.2.2. Constraints
#### 2.2.2.1. Flows balances in the hub
In each hub $h$ and each carrier $c$, the total inlet flow going to the processes bloc is either coming from the hub's import port or the return port:

$$\mathcal{F}^{in}_{c,h,t} = \mathcal{F}^{imp}_{c,h,t} + \mathcal{F}^{ret}_{c,h,t} \quad , \forall h \in N_h, \forall c \in N^{in}_h, \forall t \in N_t \tag{3}$$

Where $N^{in}_h$ is the set of carriers consumed by at least one process in hub $h$.

This same total inlet flow is the sum of all input flows of the processes for the same carrier:

$$\mathcal{F}^{in}_{c,h,t} = \sum_{p \in N^{in}_{c,h}} \mathcal{F}^{in}_{c,h,p,t} \quad , \forall h \in N_h, \forall c \in N^{in}_h, \forall t \in N_t \tag{4}$$

Where $N^{in}_{c,h}$ is the set of processes that have an inlet of carrier $c$ in hub $h$ and $\mathcal{F}^{in}_{c,h,p,t}$ is the inlet flow of carrier $c$ in process $p$ that belongs to hub $h$.

Similarly, the total output flow going out from the processes' bloc is equal to the sum of all output flows of the processes for the same carrier:

$$\mathcal{F}^{out}_{c,h,t} = \sum_{p \in N^{out}_{c,h}} \mathcal{F}^{out}_{c,h,p,t} \quad , \forall h \in N_h, \forall c \in N^{out}_h, \forall t \in N_t \tag{5}$$

Where $N^{out}_{c,h}$ is the set of processes with an outlet of carrier $c$ in hub $h$. $N^{out}_h$ is the set of carriers produced by at least one process in hub $h$.

The total output flow of carrier $c$ is added to the flow discharged from the storage component and the flow extracted from the network port. The resulting flow could be to meet the load and/or stored in the storage component and/or injected to the network and/or exported to outside the boundaries of the studied network:

$$\mathcal{F}^{out}_{c,h,t} + \mathcal{F}^{net,out}_{c,h,t} + \mathcal{F}^{sto,dis}_{c,h,t} = \mathcal{F}^{load}_{c,h,t} + \mathcal{F}^{net,in}_{c,h,t} + \mathcal{F}^{sto,cha}_{c,h,t} + \mathcal{F}^{exp}_{c,h,t} \quad , \forall h \in N_h, \forall c \in N^{out}_h \cup N^{net}_h, \forall t \in N_t \tag{6}$$

Where $N^{net}_h$ is the set of carriers for which a network port exists in hub $h$. For a carrier $c$ in a hub $h$, $\mathcal{F}^{net,in}_{c,h,t}$ and $\mathcal{F}^{net,out}_{c,h,t}$ are respectively the flows injected and extracted from the network. $\mathcal{F}^{sto,cha}_{c,h,t}$ and $\mathcal{F}^{sto,dis}_{c,h,t}$ are respectively the flows charged and discharged from the storage unit.

#### 2.2.2.2. Processes model
A generic set of constraints is used to model the processes independently of the number and nature of their inputs and outputs. The first set of constraints is used to calculate the total input flow (all carriers combined) of each process $p$:

$$\mathcal{F}^{in}_{h,p,t} = \sum_{i \in N^{in}_{h,p}} \mathcal{F}^{in}_{i,h,p,t} \quad , \forall h \in N_h, \forall p \in N^{pro}_h, \forall t \in N_t \tag{7}$$

Where $N^{pro}_h$ is the set of processes in the hub $h$ and $N^{in}_{h,p}$ is the set of inlets of the process $p$ in the hub $h$.

Similarly, the total output flow (all carriers combined) is modeled by the following set of constraints:

$$\mathcal{F}^{out}_{h,p,t} = \sum_{o \in N^{out}_{h,p}} \mathcal{F}^{out}_{o,h,p,t} \quad , \forall h \in N_h, \forall p \in N^{pro}_h, \forall t \in N_t \tag{8}$$

Where $N^{out}_{h,p}$ is the set of outlets of the process $p$ in the hub $h$. A global efficiency $\eta_g$ of the process is defined as the ratio between the total output flow and the total input flow of the process, which gives the following constraint:

$$\mathcal{F}^{out}_{h,p,t} = \eta_g \mathcal{F}^{in}_{h,p,t} \quad , \forall h \in N_h, \forall p \in N^{pro}_h, \forall t \in N_t \tag{9}$$

In this paper, a mean and constant global efficiency is considered to keep the model linear.

Furthermore, each input flow is a fraction of the total process input flow:

$$\mathcal{F}^{in}_{i,h,p,t} = f_i \mathcal{F}^{in}_{h,p,t} \quad , \forall h \in N_h, \forall p \in N^{pro}_h, \forall i \in N^{in}_{h,p}, \forall t \in N_t \tag{10}$$

A similar set of constraints is defined for the process outlets:

$$\mathcal{F}^{out}_{o,h,p,t} = f_o \mathcal{F}^{out}_{h,p,t} \quad , \forall h \in N_h, \forall p \in N^{pro}_h, \forall o \in N^{out}_{h,p}, \forall t \in N_t \tag{11}$$

Inputs and outputs fractions are also considered to be constant to keep the model linear.

2.2.2.3. Networks' model

The flows transmitted between the hubs through the networks are subject to losses. When a network port of a carrier $c$ to one or many hubs, a constraint is added to model the flows' balance in the network and transmission losses. In this paper, a simple model is adopted given by equation (12):

$$\sum_{h \in N_c^{net}} \mathcal{F}_{c,h,t}^{net,out} = (1 + f_c^{loss}) \sum_{h \in N_c^{net}} \mathcal{F}_{c,h,t}^{net,in}, \forall c \in N^{net}, \forall t \in N_t \qquad (12)$$

Where $N_c^{net}$ is the set of hubs in which a network port for carrier $c$ exist, $N^{net}$ is the set of carriers that are transmitted by a network and $f_c^{loss}$ is the network loss fraction for carrier $c$.

# 3. Case study

Geographical limitations are one of the most important factors justifying the use of multiple hubs. It is possible to identify the location of the hubs based on two main parameters influenced by geographical conditions: First, the availability of resources, and second, the demand from different carriers. To make it clear, a case study is presented here.

Accordingly, there are three different zones, and within each of them, a hub is defined, taking into consideration the characteristics of the zone:

**Hub n° 1: renewable energy hub**

This hub represents a zone near the sea. In this zone, there is sufficient space available for the placement of renewable energy power production sources (such as photovoltaic thermal hybrid solar collectors (PVT) and wind turbines (WT)) and non-compact storage systems (such as compressed air energy storage systems (CAES)). In addition, since the hub is near the sea, it is possible to use offshore or onshore wind turbines. Moreover, it provides the necessary resources (water) for desalination systems. The zone's location is far from the residential area, which may restrict access to the natural gas network. Furthermore, there is no demand for any of the carriers in this region. The details concerning the systems and the hub can be found in Table 1 and Table 2.

**Hub n° 2: residential hub**

This hub corresponds to a residential area. Due to the limited space available in this zone and the proximity of the habitats, renewable energies or systems that require large spaces may not be used. Instead, more compact and silent systems should be considered. In Table 1 and Table 2, all the processes and information concerning this hub are provided.

**Hub n° 3: industrial hub**

The industrial hub is mainly composed of offices and industries. It can provide thermal and electrical load along with potable water. Two types of thermal loads are added to the hub: industrial load, which is a high-temperature load, and offices load, which is a mid-temperature load. In addition, both the industrial sector and the offices will consume electricity. In this hub, natural gas and waste can be imported. Detailed information concerning this hub can be found in Table 1 and Table 2.

Now, after defining different hubs with their systems, a more precise identification of each system is needed. Thus, the contribution of each system is identified. Based on the developed model, for each system, the following parameters are needed: A coefficient that identifies the ratio between the total input and total output (defined in the first column of Table 1 in front of each system), a ratio that identifies the percentage of each input/output compared to total input/output (defined in the second and third columns of Table 1 in front of each input/output), a limit value for the input and/or output carrier (defined in the fourth column of Table 1).

In Table 1, since solar radiation and wind speed are both functions of time, the limit values will also vary with time. In Fig. 2, the average hourly solar radiation, wind speed, and outdoor temperature are taken from [23] are presented.

For the WT output power, the defined coefficient in Table 1 will be calculated using the WT power production equation: $P_{Electricity}^{Out\ WT} = \frac{\pi}{2} r^2 \rho \eta v^3 = Coefficient \times v^3$. By using this coefficient and data of Fig. 2, the WT output power can be calculated. In this equation, $r$ is the blade radius, $\rho$ is the air density, $\eta$ is the WT efficiency, and $v$ is the wind speed. The cut-in and cut-off speeds are $3\ m/s$ and $25\ m/s$, respectively. For the PVT system, the solar radiation in W will also be calculated as $PVT\ area\ (m^2) \times Average\ solar\ radiation\ (W/m^2)$ using Fig. 2. Storage systems and ports added to each hub are presented in Table 2.

All hubs have a mid-temperature export port with a price equal to zero. This allows excess heat to be rejected. Indeed, in some cases, electricity production (using CHP, fuel cell, or PVT) will generate heat that cannot be consumed or stored. If the choice and size of the processes are not properly established, this heat rejection can have a high value. Moreover, as previously discussed, the waste entering the industrial hub has a constant flow and must be incinerated in the hub.

*Table 1. Systems identification*

| Processes ($Total\ outlet/Total\ inlet$) | Inputs (ratios) | Outputs (ratios) | Limits |
|---|---|---|---|
| **Hub n° 1: renewable energy hub** | | | |
| Fuel cell (0.75) | $H_2^{High\ pressure}\ W$ (1) | Electricity $W$ (2/3) <br> High-temperature heat $W$ (1/3) | $Power_{Electricity}^{Max\ out} = 200\ kW$ |
| Electrolyzer (0.9) | Electricity $W$ (1) | $H_2^{Low\ pressure}$ (1) | $Power_{H_2}^{Max\ out} = 400\ kW$ |
| H2 compressor (0.9) | Electricity $W$ (0.1) <br> $H_2^{Low\ pressure}\ W$ (0.9) | $H_2^{High\ pressure}\ W$ (1) | $Power_{H_2}^{Max\ out} = 400\ kW$ |
| Desalination (650 $kg/kW$) | High-temperature heat $W$ (1) | Potable water $kg$ (1) | $Production_{Potable\ water}^{Max\ out} = 10^5\ kg/s$ |
| PVT (0.5) <br> Area = $20\ m^2$ | Solar radiation $W$ (1) | Electricity $W$ (1/5) <br> Mid-temperature heat $W$ (4/5) | $P_{Solar\ radiation}^{Max\ in} = Variable\ by\ time\ W$ |
| WT (5 $kW\ s^3/m^3$) | Wind $m/s$ (1) | Electricity $W$ (1) | $P_{Electricity}^{Max\ out} = Variable\ by\ time\ W$ |
| **Hub n° 2: residential hub** | | | |
| Gas boiler (0.9) | Natural gas $W$ (1) | Mid-temperature heat $W$ (1) | $P_{Heat}^{Max\ out} = 2 \times 10^3\ kW$ |
| Mid-temperature HP (4) | Electricity $W$ (1) | Mid-temperature heat $W$ (1) | $P_{Electricity}^{Max\ in} = 5 \times 10^2\ kW$ |
| **Hub n° 3: industrial hub** | | | |
| Waste furnace (0.85) | Waste $W$ (1) | High-temperature heat $W$ (1) | $P_{Heat}^{Max\ out} = 5 \times 10^3\ kW$ |
| High-temperature HP (1) | Electricity $W$ (0.25) <br> Mid-temperature heat $W$ (0.75) | High-temperature heat $W$ (1) | $P_{Electricity}^{Max\ in} = 7 \times 10^2\ kW$ |
| CHP (0.6) | Natural gas $W$ (1) | Electricity $W$ (7/12) <br> Mid-temperature heat $W$ (5/12) | $P_{Electricity}^{Max\ out} = 7.2 \times 10^3\ kW$ |
| Heat exchanger (0.9) | High-temperature heat $W$ (1) | Mid-temperature heat $W$ (1) | $P_{Heat}^{Max\ out} = 1.5 \times 10^3\ kW$ |

In Table 2, one of the most significant parameters is the load. To accurately identify the operating mode of the systems, load profiles are the most important parameters. Now, different hourly load profiles for the first week of January, from Monday to Sunday, are presented.

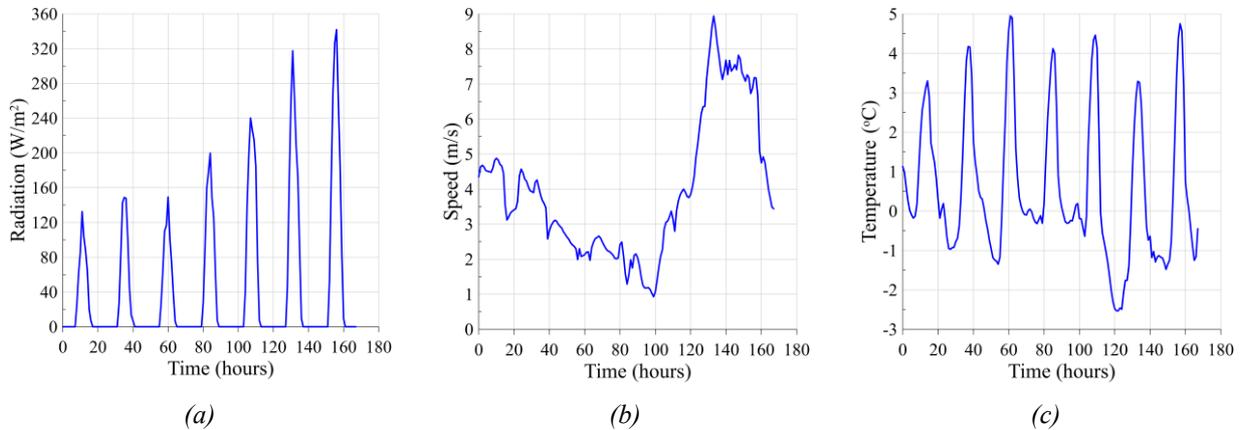

*Fig. 2. Metrological data* [24]: *a) hourly solar radiation b) hourly wind speed, c) Outdoor temperature*

*Table 2. Hubs' identification*

| Storages (Efficiency %) | Import ports (Price €/kW) | Return ports | Network ports (Maximum flow) | Export ports (Price €/kW) | Loads |
|---|---|---|---|---|---|
| **Hub n° 1: renewable energy hub** | | | | | |
| Water tank (100) | Solar radiation (0) | Electricity | Electricity ($3 \times 10^3\ W$) | Mid-temperature heat (0) | |
| H2 tank (97) | Wind (0) | High-temperature heat | High-temperature heat ($2 \times 10^2\ W$) | | |
| CAES (65) | | $H_2^{Low\ pressure}$ | Mid-temperature heat ($2.5 \times 10^3\ W$) | | |
| | | $H_2^{High\ pressure}$ | Potable water ($1,4 \times 10^5\ kg/s$) | | |
| **Hub n° 2: residential hub** | | | | | |
| Thermal storage (80) | Natural gas (10) | Electricity | Electricity ($5.5 \times 10^3\ W$) | Mid-temperature heat (0) | Electricity |
| Battery (85) | | | Mid-temperature heat ($5 \times 10^3\ W$) | | Mid-temperature heat |
| | | | Potable water ($1,2 \times 10^5\ kg/s$) | | Potable water |
| **Hub n° 3: industrial hub** | | | | | |
| | Natural gas (10) | Electricity | Electricity ($6.3 \times 10^3\ W$) | Mid-temperature heat (0) | Electricity |
| | Waste (0) | High-temperature heat | High-temperature heat ($6 \times 10^2\ W$) | | High-temperature heat |
| | | Mid-temperature heat | Mid-temperature heat ($5 \times 10^3\ W$) | | Mid-temperature heat |
| | | | Potable water ($2,2 \times 10^4\ kg/s$) | | Potable water |

In Fig. 3, all the electrical and thermal loads are illustrated. The residential electrical load is shown in Fig. 3.a. It should be noted that the residential population is 5000 persons. Furthermore, for the electrical load of the offices, it is assumed that the electrical load per capita will be 30% of the residential one. The offices' population is equal to 2500 persons.

For the industrial electrical load, 20 industries with similar electrical load profiles in Fig. 2.b are considered.

For the potable water load in the residential hub, the data of [24] is used. Reference [24] gives the daily potable water load per capita and also the distribution of this potable water utilization per hour. The distribution of potable water utilization as a fraction of total daily consumption is presented in Fig. 2.c.

Fig. 3.c, which provides the hourly break-up of the total consumed potable water, and daily water consumption of 352 L/day/capita will provide the hourly water consumption per capita. Also, to calculate the fraction of the offices' water consumption, the fraction of the water consumption of the residential one which does not exist in the offices (e.g., shower and clothes wash) will be eliminated, and the rest will be considered as the fraction of the water consumption of the offices. The results will be 59% of the total residential potable water consumption (352 L/day/capita) with the daily fractions presented in Fig. 3.d.

For the residential and offices' thermal consumption, an indoor set-point temperature is assumed and combined with heat transfer coefficients and outdoor temperature (Fig. 1.c) to calculate the heat load. For the residential, it assumed that from 8:00 to 22:00, the inside temperature is equal to 20℃ and outside these hours is equal to 19℃. For the offices, the inside temperature is kept at 20℃ between 7:00 and 19:00, and for the weekend or outside these hours, it is assumed to be equal to 15℃. By using the thermal coefficient of 0.06 $kW/m^2K$ [25] and indoor and outdoor temperatures, the heat load profiles for the residential sector and offices can be extracted as Fig. 2.e. It should be noted that the residential and offices areas are $10^5\ m^2$ and $2.5 \times 10^3\ m^2$, respectively.

The thermal load for the industrial sector is also shown in Fig. 3.f. It should be noted that 32 industries with the thermal load profile of Fig. 2.f exist. The industrial thermal load needs to be satisfied by high-temperature heat.

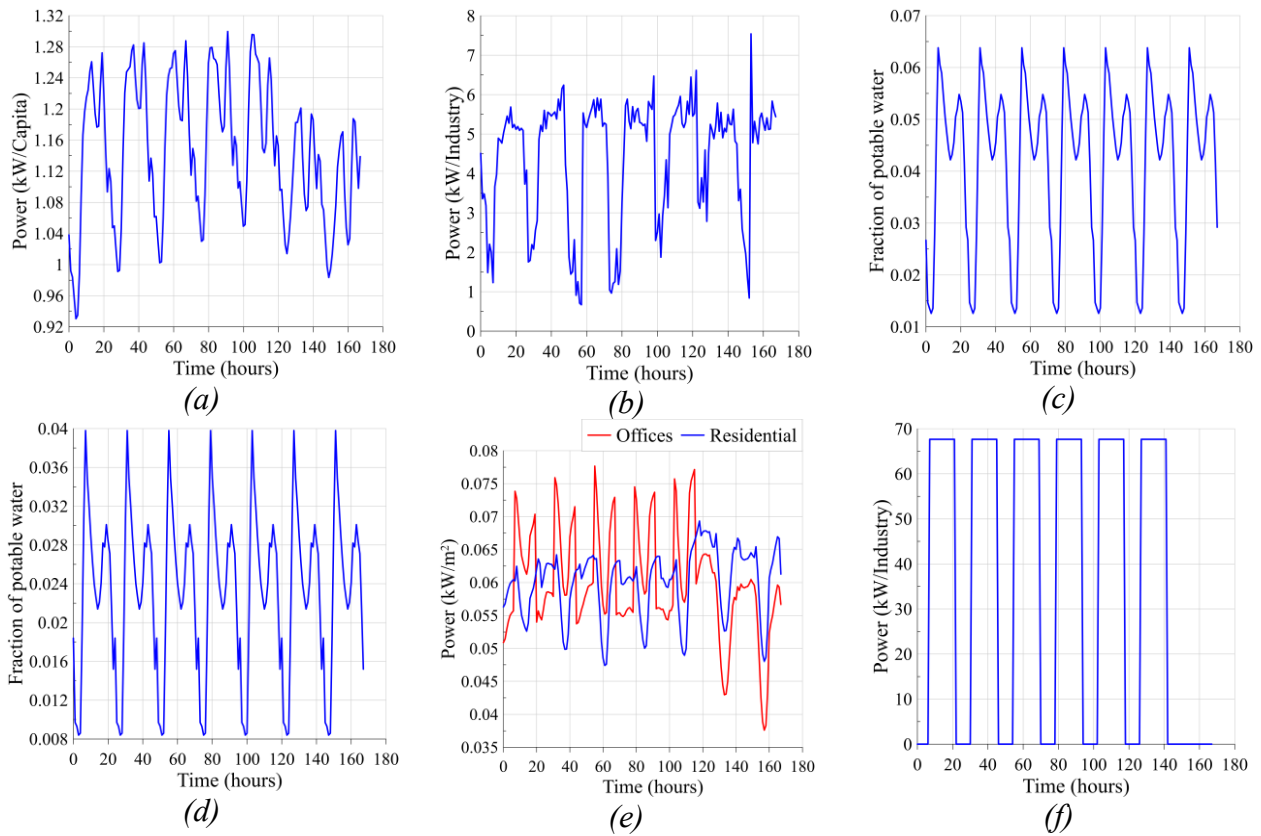

Fig. 3. Load data: a) Residential electrical load per capita [26], b) Industrial electrical load per industry [27], c) Residential potable water consumption as a fraction of total daily consumption [24], d) Offices' potable water consumption as a fraction of total daily consumption, e) Heating thermal load for residential and offices, f) Thermal load per industry [28]

## 4. Results and discussions

Optimal scheduling applied to the case study gives the flows and systems utilization to minimize the operational energy cost.

The systems and storages behaviors of different hubs are shown in Fig. 4- Fig. 6 for the three hubs.

The first point that should be noted here (from Fig. 4.c and Fig. 4.d), is that the WT and PVT are working at their maximum capacities. There is no curtailment observed for these systems. This is possible thanks to the use of storage systems and the synergy between the different processes. For the furnace system, since the waste flow is constant and must be used, the furnace production is constant.

Other systems that can help produce electricity are CHP (Fig. 6.c) and fuel cell (Fig. 4.a). It seems even though the CHP is using natural gas, it is preferred to the fuel cell. This behavior is explained by the poor global efficiency of the hydrogen cycle. Indeed, the hydrogen should be produced using the electrolyzer and then pressurized by a compressor, and both systems consume electricity that should be produced by WT, PVT, or CHP. This justifies that fuel cell is not the better option for electrical and thermal production than the CHP system. It can also be seen in Fig. 4.c, Fig. 4.d, and Fig. 5.c that WT and PVT power production increase causes CHP production to decrease.

The rest of the high-temperature thermal energy demand is produced by the high-temperature HP located in hub 3 (Fig. 6.b). This system rises the temperature of a part of the mid-temperature thermal energy produced either by the CHP in the same hub, the mid-temperature HP, the gas boiler in hub 2, or the PVT in hub 1. Since there is a furnace (Fig. 6.a) that works all the time, and there are not any high-temperature storage units, it is logical that a part of the produced high-temperature thermal energy will be used to satisfy mid-temperature demands through the heat exchanger (Fig. 6.d). The remaining part of the mid-temperature thermal energy demand is produced by the mid-temperature HP (Fig. 5.b) and gas boiler (Fig. 5.a). Since the input of the mid-temperature HP is electricity and the input of the gas boiler is the natural gas, a trade-off between these systems is made by the optimizer.

Fig. 4 and Fig. 5 show that all the storages are empty at the end of the week. This is logical since the optimization is done on for one week. Furthermore, Fig. 4.g and Fig. 5.d show that the battery, due to its higher efficiency (Table 2), is more adapted for the short-term utilization, while the CAES is preferred to be used for long-term usage with fewer charge and discharge cycles. Indeed, CAES and hydrogen storage are discharged mainly during the period between 60h and 100h when there is no wind production and between 160h and 180h when wind production decreased drastically. Outside these periods, these systems are either charging or kept fully charged.

The total quantities transferred from each carrier for each port are shown in Table 3. It indicates that some of the thermal export ports are not used. This is the case in renewable energy and residential hubs. Also, the exported thermal energy from the industrial hub is not too high, demonstrating that the optimization's operational strategy reduces energy waste. For the network ports, the negative sign indicates carrier import from the network and the positive sign means carrier export to the network. For example, since the potable water will be produced only by the first hub, the other hubs will import it. Furthermore, the high-temperature thermal energy produced at the renewable energy hub is not enough to satisfy the input energy for the desalination system to produce the required potable water. Thus, importing high-temperature thermal energy into the renewable energy hub will be unavoidable. In addition, the return port shows the carriers recycled inside the hubs. For example, in the industrial hub, mid-temperature heat is used by the high-temperature HP to produce high-temperature heat. Also, high-temperature heat is used by the heat exchanger that downgrades the temperature to fulfill the space heating demand of the offices.

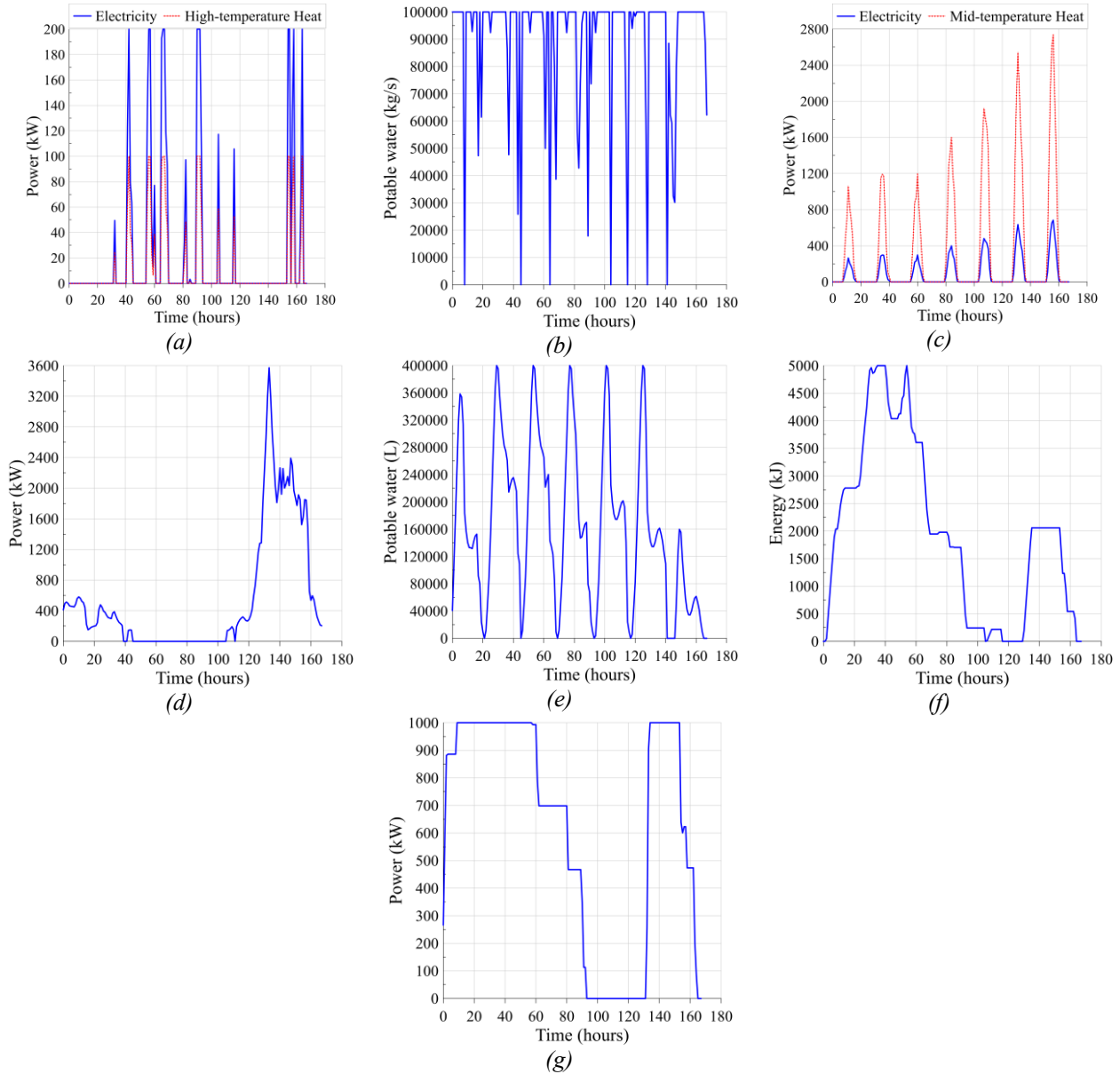

*Fig. 4. Results of hub n°1 (renewable energy): a) Fuel cell production, b) Desalination system production, c) PVT production, d) WT production, e) Water tank stored water, f) Hydrogen tank stored energy g) CAES stored energy*

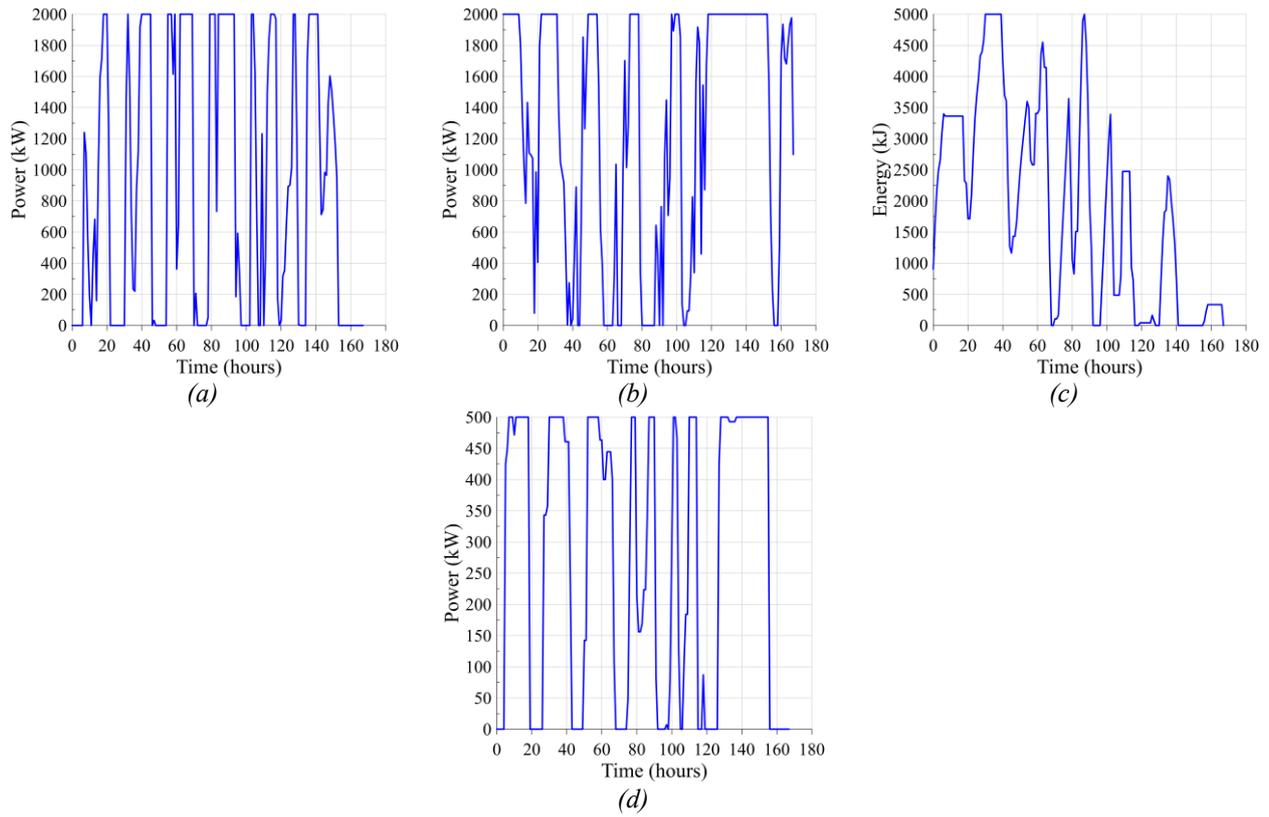

Fig. 5. Results of hub n°2 (residential): a) Gas boiler production, b) Mid-temperature HP, c) Mid-temperature thermal storage, d) Battery

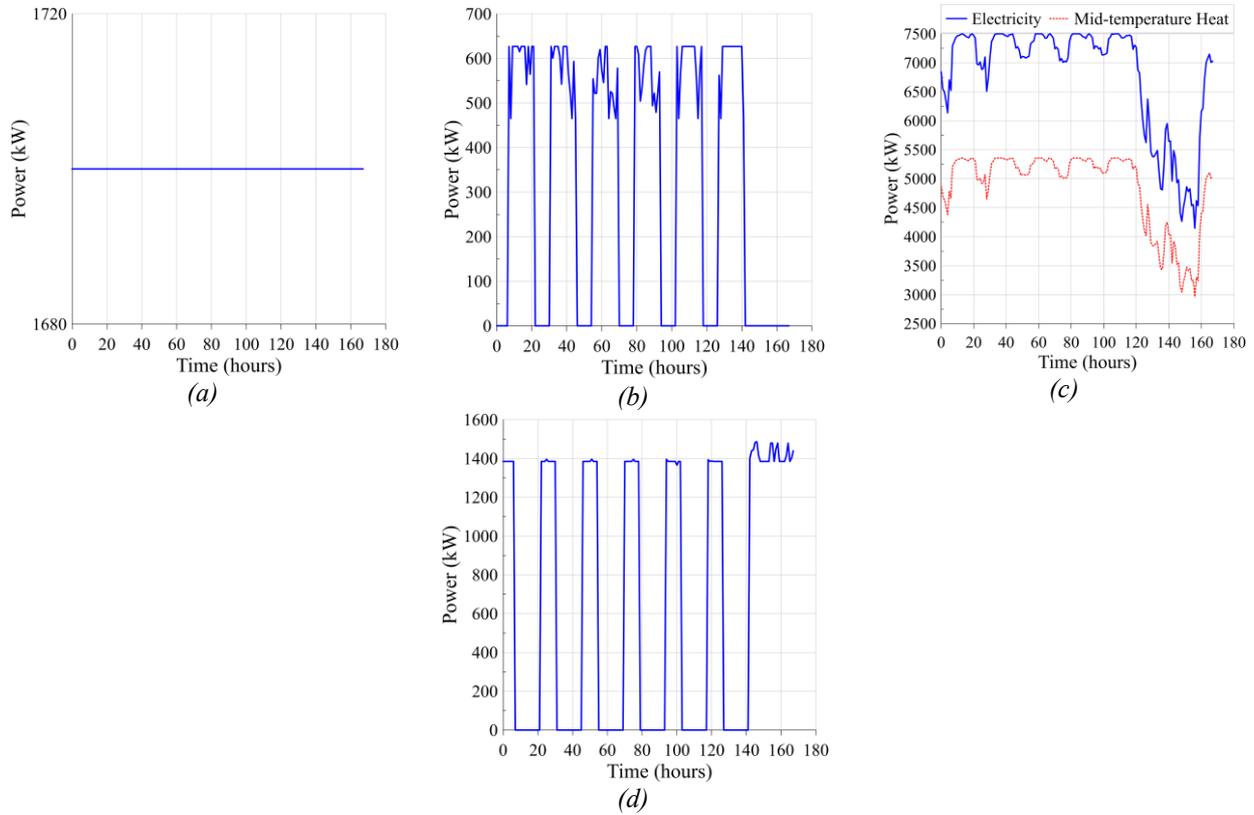

Fig. 6. Results of hub n°3 (industrial): a) Waste furnace production, b) High-temperature HP, c) CHP, d) Heat exchanger

*Table 3. Total port utilization of different hubs*

| Import ports (Total $kWh$) | Return ports (Total $kWh$) | Network ports (Total $kWh$) | Export ports (Total $kWh$) |
|---|---|---|---|
| **Hub n° 1: renewable energy hub** | | | |
| Solar radiation ($1.505642 \times 10^5$) | Electricity ($2.26236 \times 10^4$) | Electricity ($9.9005 \times 10^4$) | Mid-temperature heat (0) |
| Wind ($1.86575 \times 10^4$) | High-temperature heat ($2.31793 \times 10^4$) | High-temperature heat ($-2.11488 \times 10^4$) | |
| | $H_2^{Low\ pressure}$ ($8.632 \times 10^3$) | Mid-temperature heat ($6.02257 \times 10^4$) | |
| | $H_2^{High\ pressure}$ ($8.1219 \times 10^3$) | Potable water ($1.50665 \times 10^7\ kg$) | |
| **Hub n° 2: residential hub** | | | |
| Natural gas ($1.667778 \times 10^5$) | Electricity ($5.36395 \times 10^4$) | Electricity ($-1.018365 \times 10^6$) | Mid-temperature heat (0) |
| | | Mid-temperature heat ($-6.55574 \times 10^5$) | |
| | | Potable water ($-1.232 \times 10^7\ kg$) | |
| **Hub n° 3: industrial hub** | | | |
| Natural gas ($3.2650756 \times 10^6$) | Electricity ($1.3132 \times 10^4$) | Electricity ($9.703532 \times 10^5$) | Mid-temperature heat ($5.6852 \times 10^3$) |
| Waste ($3.36 \times 10^5$) | High-temperature heat ($1.21043 \times 10^5$) | High-temperature heat ($2.22271 \times 10^4$) | |
| | Mid-temperature heat ($3.93961 \times 10^4$) | Mid-temperature heat ($6.28127 \times 10^5$) | |
| | | Potable water ($-2.02908 \times 10^6\ kg$) | |

## Conclusion

In this paper, a novel formalism is presented to model networked multi-carrier hubs. The formalism is illustrated mathematically using a linear optimization problem that optimizes the operation cost of existing hubs. The model formulation could be easily extended to other objective functions and nonlinear programming approaches while keeping the same architecture. The formalism is implemented in a framework allowing to model a network of multi-carrier hubs by simply adding the processes, storage devices, loads, and different exchange ports. The framework is used to model and find the optimal operation strategy of a theoretical case study consisting of three hubs involving different carriers: Electricity, gas, heat, hydrogen, water, etc.

Further developments are in progress to add new functionalities to the framework, such as optimal processes selection, custom objective functions, ramping limitations, stop and start scheduling, nonlinear behavior of the processes, etc.

## Nomenclature

| | | | | | |
|---|---|---|---|---|---|
| $c$ | cost, € | $\dot{m}$ | mass flow rate, kg/s | **Subscripts and superscripts** | |
| $p$ | price, € | $N$ | set | $in$ input | $h$ hub |
| $f$ | fraction, - | $\eta$ | efficiency | $out$ output | $c$ carrier |
| $\mathcal{F}$ | carrier flow, J/s | | | $pro$ process | $p$ process |
| | | | | | $t$ time |